\documentclass[11pt,reqno]{amsart}
\usepackage{amssymb}
\input amssym.def
\usepackage{amsmath, amsfonts}
\usepackage{amssymb}
\usepackage{amscd, mathtools}
\usepackage[mathscr]{eucal}
\usepackage[usenames]{color}
\usepackage{fancyhdr}
\usepackage{textpos}

\newfont{\cyrr}{wncyr10}

\newcommand{\Z}{{\mathbb Z}}
\newcommand{\Q}{{\mathbb Q}}
\newcommand{\R}{{\mathbb R}}
\newcommand{\C}{{\mathbb C}}

\newcommand{\K}{{\mathbf K}}
\newcommand{\T}{{\rm Tr}}

\newcommand{\F}{{\mathbf F}}
\renewcommand{\L}{{\mathbf L}}

\newcommand{\thmref}[1]{Theorem~\ref{#1}}
\newtheorem{thm}{Theorem}

\newtheorem{lem}[thm]{Lemma}
\newtheorem{cor}[thm]{Corollary}
\newtheorem{prop}[thm]{Proposition}
\newtheorem{rmk}[thm]{Remark} 
\newtheorem{defn}[thm]{Definition}
\newtheorem{conj}[thm]{Conjecture}

\newcommand{\propref}[1]{Proposition~\ref{#1}}
\newcommand{\lemref}[1]{Lemma~\ref{#1}}
\newcommand{\defref}[1]{Definition~\ref{#1}}

\begin{document}

\fancypagestyle{plain}{%
	\fancyhead[R]{\fbox{to appear in journal xx}}
	\renewcommand{\headrulewidth}{1pt}
}

\title[Linear independence of Dirichlet L values]{On linear 
independence of Dirichlet $L$ values}

\author{Sanoli Gun, Neelam Kandhil and Patrice Philippon}

\address{Sanoli Gun and Neelam Kandhil \\ \newline
The Institute of Mathematical Sciences, A CI of Homi Bhabha National Institute, 
CIT Campus, Taramani, Chennai 600 113, India.}

\address{Patrice Philippon \\ \newline
Institut de Math\'ematiques de Jussieu, UMR 7586 du CNRS, 
4 place Jussieu, 75252 Paris Cedex 05, France.}

\email{sanoli@imsc.res.in} 
\email{neelam@imsc.res.in} 
\email{patrice.philippon@upmc.fr}

\subjclass[2010]{11J72, 11R18, 11M06}
\keywords{Dirichlet $L$-functions, cotangent values, cyclotomic field, Okada's theorem}

\begin{abstract} 
The study of linear independence of $L(k, \chi)$ for a fixed integer $k>1$ 
and varying $\chi$  depends critically on the parity of $k$ vis-\`a-vis $\chi$. 
This has been investigated by a number of authors for Dirichlet characters 
$\chi$ of a fixed modulus and having the same parity as $k$.
The focal point of this article is to extend this investigation to 
 families of Dirichlet characters   modulo distinct pairwise co-prime natural numbers. 
 The interplay between the resulting ambient number fields brings in
 new technical issues and complications hitherto absent in the context
 of a fixed modulus (consequently a single number field lurking in the background).
 This entails a very careful and hands-on  dealing with the arithmetic of compositum
 of number fields  which we undertake in this work. Our results extend earlier works 
of the first author with Murty-Rath as well as  works of Okada, Murty-Saradha and Hamahata.
\end{abstract}

\maketitle

\section{\large{Introduction and statements of Theorems}}
For a Dirichlet character $\chi$ modulo $q >1$ and $s \in \C$, consider the
Dirichlet $L$-function
$$
L(s, \chi) = \sum_{n=1}^{\infty} \frac{\chi(n)}{n^s} 
$$
which converges absolutely for $\Re(s) > 1$. Further, it is holomorphic
in this region. The study of irrationality of $L(k, \chi)$ for natural numbers $k>1$ 
has an intriguing history starting from the work of Euler. He found a
closed formula for $L(k, \chi_0)$ when $\chi_0$ is the trivial character 
modulo $q \ge 1$ and $k$ is even. 
When $k$ is odd, it follows from the work of Ball and 
Rivoal \cite{BR} (see also \cite{FSZ, LY}) that there 
are infinitely many irrational numbers as $k>1$ varies over odd natural numbers. 
When $\chi$ is the non-trivial character modulo $4$, similar results were
established by Rivoal and Zudilin \cite{RZ}. For arbitrary non-trivial 
character $\chi$ modulo $q$, infinitude of irrationality of $L(k, \chi)$  
when $\chi(-1) = (-1)^{k+1}$ follows from the recent work of Fischler \cite{SF}.
When $\chi(-1) = (-1)^{k}$,  infinitude of irrationality of $L(k, \chi)$  is well
known (see \cite[Ch. VII, \S2]{JN}).

In this article, we study linear independence of
$L(k, \chi)$ when $k$ is fixed and $\chi$ varies over Dirichlet characters
modulo pairwise co-prime natural numbers.
If we fix $k>1$ and vary $\chi$ modulo a natural number $q>2$, the
question of linear independence of $L(k, \chi)$'s over $\Q$ was first 
investigated by Okada \cite{TO}.  As noted by Murty-Saradha \cite{MS1},
this result can be extended over number fields which are 
disjoint to the $q$th cyclotomic field. To proceed further, we need
to introduce some notations. For a natural number $k>1$, let us 
set 
$$
X_{q, k} = \{ L(k,\chi) ~|~ \chi \bmod{q},  ~\chi \ne \chi_0 \}.
$$
We can write $X_{q, k} = X_{q, k, e} \cup X_{q, k, o}$, where
\begin{eqnarray}\label{eos}
	X_{q, k, e} 
	&=& \{ L(k,\chi) ~|~ \chi \bmod{q}, ~ \chi(-1) = 1, \chi \ne \chi_0 \} \nonumber \\
	\phantom{mm}\text{and}\phantom{mm}
	X_{q, k, o} 
	&=& \{ L(k,\chi) ~|~ \chi \bmod{q},  ~\chi(-1) = -1 \}.
\end{eqnarray}
Further, for any non zero natural number $q$, let $\zeta_q$ denotes a
primitive $q$th root of unity. 
In this set-up, Okada \cite{TO} (see also Murty-Saradha \cite{MS1}) proved 
the following theorems.
\begin{thm}{\rm(Okada \cite{TO})}
Let $k \geq 1, q>2$ be natural numbers and $\K$ be a number field with 
$\K(\zeta_{\varphi(q)}) \cap \Q(\zeta_q) = \Q$.  Then the numbers in the set 
$X_{q, 2k+1, o}$ are linearly independent over $\K(\zeta_{\varphi(q)})$. 
\end{thm}

\begin{thm}{\rm(Okada \cite{TO})}
Let $k \geq 1, q>2$ be natural numbers and $\K$ be a number field with 
$\K(\zeta_{\varphi(q)}) \cap \Q(\zeta_q) = \Q$. Then the numbers in the set 
$\{\zeta{(2k)}\} \cup X_{q, 2k, e}$
are linearly independent over $\K(\zeta_{\varphi(q)})$. 
\end{thm}

From now on, for an integer $q>2$, we shall denote the maximal
real subfield of $\Q(\zeta_{q})$ by $\Q(\zeta_{q})^{+} $. 
In this set-up, we have the following theorems.
\begin{thm}\label{thm1}
For $1 \le j \le \ell$, let $q_j>2$  be pairwise co-prime natural numbers and 
$k \geq 1$ be an integer. If $\K$ is a number field such that 
$\K(\zeta_{\varphi(q_1) \cdots \varphi(q_{\ell})}) \cap  \Q(\zeta_{q_1\cdots q_{\ell} })^{+} = \Q$,
then the numbers in the set 
$$
X_{q_1, 2k+1, o} \cup \cdots \cup X_{q_\ell, 2k+1, o}
$$
are linearly independent over $\K(\zeta_{\varphi(q_1) \cdots \varphi(q_{\ell})})$. 
\end{thm}
\begin{thm}\label{thm2}
For $1 \le j \le \ell$, let $q_j>2$  be pairwise co-prime natural numbers and $k \geq 1$ be an integer.
If $\K$ is a number field such that 
$\K(\zeta_{\varphi(q_1) \cdots \varphi(q_{\ell})}) \cap \Q(\zeta_{q_1\cdots q_{\ell} })^{+} = \Q$, 
then the numbers in the set 
$$
\{\zeta{(2k)}\} \cup X_{q_1,2k, e} \cup \cdots \cup X_{q_\ell, 2k, e}
$$
are linearly independent over $\K(\zeta_{\varphi(q_1) \cdots \varphi(q_{\ell})})$. 
\end{thm}

If one replaces Dirichlet characters by arbitrary 
periodic arithmetic functions $f$ with period $q>1$ and considers the
associated $L$-function
$$
L(s,f) = \displaystyle \sum_{n=1}^{\infty} \frac{f(n)}{n^s}
$$
for $s \in \C$ with $\Re(s) >1$,  then non-vanishing of
$L(k, f)$ is intricately related to a conjecture of Chowla-Milnor \cite{GMR}.
For describing this conjecture, we need to introduce Hurwitz zeta function
$$
\zeta(s,x) = \displaystyle \sum_{n=0}^{\infty} \frac{1}{(n + x)^s}~,
$$
where $x$ is a real number with $0 < x \leqslant 1$ and $s$ is a 
complex number with $\Re(s) > 1$. It is easy to see that when 
$k >1$, we have
\begin{equation}\label{eq1}
L(k,f) = q^{-k} \displaystyle \sum_{a=1}^{q} f(a)\,\zeta(k, a/q)
\end{equation}
and hence in particular, when $\chi_0$ is the trivial character modulo $q\ge2$,
\begin{equation}\label{eq2}
\zeta(k) \prod_{\stackrel{p|q}{\scriptscriptstyle p\text{ prime}}}  (1-p^{-k})  = 	
L(k, \chi_0)=
 q^{-k} \sum_{1 \leq a < q \atop (a,q)=1 } \zeta(k, a/q).
\end{equation}
For example, $\zeta(k,1/2) = (2^k-1)\zeta(k) \not= 0$, for all $k>1$.

\begin{rmk}
Since $\zeta(s,x)$ extends analytically
to the entire complex plane, apart from $s=1$, where it has a
simple pole with residue $1$, we have by \eqref{eq1} that
$L(s,f)$, for a periodic function $f$ modulo $q$, 
extends meromorphically to the complex plane with a
possible simple pole at $s=1$ with residue $q^{-1}\sum_{a=1}^{q} f(a)$.
Thus when $f =\chi$, a non-trivial Dirichlet character modulo $q$, the 
number $L(1, \chi)$ makes sense. See the articles \cite{BBW, GK, NR, MM1, TO}
for linear independence of such values. 
\end{rmk}

P. Chowla and S. Chowla  \cite{CC} were the first to study non-vanishing of
$L(2,f)$ for arbitrary periodic functions $f$ and made the following conjecture.
\begin{conj}{\rm (Chowla-Chowla)}
Let $p$ be any prime and $f$ be any rational valued periodic function with 
period $p$. Then $L(2,f) \neq 0$
except in the case when
$$
f(1) = f(2) = \cdots = f(p-1) = \frac{f(p)}{1-p^2}~.
$$
\end{conj}

Milnor \cite{MIL} reformulated the conjecture of Chowla-Chowla as follows;
\begin{conj}{\rm (Milnor).} 
For any integer $k >1$, the real numbers 
$$
\zeta\big(k, 1/p \big), \phantom{b} \zeta\big(k, 2/p \big), \cdots ,
\zeta\big(k,(p-1)/p \big) 
$$
are all linearly independent  over $\Q$.
\end{conj}
When $q$ is not necessarily prime, Milnor suggested the 
following generalization of the Chowla conjecture.

\begin{conj}{\rm (Chowla-Milnor)} \label{conjCM}
Let $k>1,q>2$ be integers. Then the following $\varphi(q)$ real numbers  
$$
\zeta(k, a/{q}) \phantom{bb} {\rm with}~~(a,q) = 1, ~~1 ~\leq a ~< q
$$ 
are linearly independent over $\Q$. 
\end{conj}

In relation to the Chowla-Milnor conjecture, we define 
the following linear spaces (see \cite{GMR}).

\begin{defn}
For a number field $\K\subset\C$ and integers $k>1$, $q\ge1$, 
the $\K$-vector space 
$$
V_{\K, k}(q) = \K- \text{span of} ~ \{ \zeta(k,a/q) : 1 \leq a \le q, (a,q)=1\}
$$
is defined to be the Chowla-Milnor space for $\K$ and $q$. In particular, $V_{\K,k}(1) 
= \K\zeta(k,1) = \K\zeta(k)$ and $V_{\K,k}(2) = \K\zeta(k,1/2) = \K\zeta(k)$.
\end{defn}
Conjecture \ref{conjCM} is equivalent to $\dim_{\Q}(V_{\Q,k}(q)) = \varphi(q)$ for $q>2$. 
We observe that $\sum_{d\mid q}V_{\K,k}(d) = \sum_{a=1}^q\K \zeta (k, a/q)$.

For $q\ge1$, we can write the space $V_{\K,k}(q)$ as 
$V_{\K,k}(q) = V^+_{\K,k}(q) + V^-_{\K,k}(q)$, where for $q>2$
\begin{eqnarray*}
V^\pm_{\K,k}(q) 
&=& 
\sum_{1 \le a < q/2 \atop{ (a, q)=1 }} \K\left( \zeta (k, a/q)\pm (-1)^k\zeta (k,1-a/q) \right) \\
\phantom{m}\text{and}\phantom{mmm}
V^\pm_{\K,k}(2) 
&=& \K\zeta (k, 1/2)(1\pm (-1)^k), \quad\quad
V^\pm_{\K,k}(1) = \K\zeta (k, 1)(1\pm (-1)^k).
\end{eqnarray*}

For $q=1,2$, we have $\dim_{\K}(V^{\pm}_{\K,k}(1)) = \dim_{\K}(V^{\pm}_{\K,k}(2)) = \frac{1}{2}(1\pm(-1)^k)$.

For $q>2$, it results from Okada's theorem \ref{okada} that 
$\dim_{\Q}(V^+_{\Q,k}(q)) = \varphi(q)/2$.

 Since $\dim_{\Q}(V^-_{\Q,k}(q)) \le \varphi(q)/2$, Conjecture \ref{conjCM} is equivalent to 
 $V^+_{\Q,k}(q) \cap V^-_{\Q,k}(q) = 0$ and $\dim_{\Q}(V^-_{\Q,k}(q)) = \varphi(q)/2$.
In this set-up, we have the following theorem.

\begin{thm}\label{cor1}
For $1 \le j \le \ell$, let $q_j\ge1$  be pairwise co-prime natural numbers and $k>1$ be an integer.
If $\K$ is a number field such that $\K \cap \Q(\zeta_{q_1\cdots q_{\ell} })^{+} = \Q$, then 
\begin{equation*}
\dim_{\K} \left({\sum}_{j=1}^{\ell} V^+_{\K, k}(q_j)\right)
= {\sum}_{j=1}^{\ell}\dim_{\K}(V^+_{\K, k}(q_j)) - (\ell-1)\dim_{\K}(V^+_{\K, k}(1)).
\end{equation*}
In particular,
\begin{equation*}
\dim_{\K} \left({\sum}_{j=1}^{\ell} V^+_{\K, k}(q_j)\right) = \begin{cases}
\sum_{j=1}^{\ell} \frac{\varphi(q_j)}{2} - \frac{\ell-1}{2}(1+(-1)^k) &\text{if } 2<q_1,\cdots,q_\ell,\\[2mm]
\sum_{j=2}^{\ell} \frac{\varphi(q_j)}{2} - \frac{\ell-2}{2}(1+(-1)^k) &\text{if } 2=q_1<q_2,\cdots,q_r.
\end{cases}
\end{equation*}
\end{thm}

When $\K =\Q$, Chowla-Milnor conjecture predicts that the dimension of 
$V_{\Q, k}(q)$ over $\Q$ is equal to $\varphi(q)$. Here we have the
following corollary.

\begin{cor}\label{cor3}
For $1 \le j \le \ell$, let $q_j>2$  be pairwise co-prime natural numbers and $k > 1$ be an integer.
If $\K$ is a number field such that $\K \cap \Q(\zeta_{q_1\cdots q_{\ell} })^{+} = \Q$, then 
$$
\frac12\sum_{j=1}^{\ell} \varphi(q_j) - \frac{\ell-1}2(1+(-1)^k) 
~\leq~ 
\dim_{\K} \left({\sum}_{j=1}^{\ell} V_{\K, k}(q_j)\right) 
~\leq~ 
\sum_{j=1}^{\ell} \varphi(q_j) - (\ell - 1).
$$
\end{cor}

 \begin{rmk}
 Let $k,q>1$ be integers and 
 $\{\chi_b  : b\in(\mathbb Z/q\mathbb Z)^\times \}$ be the set of
 Dirichlet characters modulo $q$.  For $a,b$ running over 
$(\mathbb Z/q\mathbb Z)^\times$, we have
\begin{equation}\label{rain}
q^kL(k,\chi_b) 
= \sum_{(a,q)=1}\chi_b(a)\zeta(k, a/q). 
\end{equation}
By the orthogonality relations satisfied by Dirichlet characters, 
the matrix $\left(\chi_b(a)\right)_{a,b}$ has inverse 
$$
\frac{1}{\varphi(q)}\left(\chi_b(a^{-1})\right)_{b,a}.
$$ 
Let $\K\subset\mathbb{C}$ be a number 
field containing the $\varphi(q)$-th roots of unity, the $\K$-vector spaces
$$
\sum_{(b,q)=1}\K L(k,\chi_b)
\phantom{m}\mbox{and}\phantom{m}
\sum_{(a,q)=1}\K \zeta (k, a/q) = V_{\K,k}(q)
$$
are equal. Furthermore, it follows from \eqref{rain} that 
\begin{equation}\label{dry}
\sum_{\chi(-1)=\pm(-1)^k}\K L(k,\chi) = V_{\K,k}^\pm(q).
\end{equation}
\end{rmk}
 
 \begin{defn}\label{productofvectorsspaces}
Let $\K$ be a number field, $V$ and $W$ be two $\K$-vector spaces in $\C$. 
We define the product $VW$ as the $\K$-span of the set of numbers $vw$ 
with $v\in V$ and $w\in W$.
\end{defn}

Following Hamahata \cite{YH}, we consider generalized Chowla-Milnor spaces.
 \begin{defn}\label{gCM}
 Let $k_1, \cdots, k_r >1$ and $q_1, \cdots, q_r \ge 1$ be integers. Set $\vec{k}= (k_1, \cdots, k_r)$
 and $\vec{q} = (q_1, \cdots, q_r)$. For a number field $\K\subset\C$,
 the generalized Chowla-Milnor space is defined by
 $$
 V_{\K, \vec{k} }( \vec{q} )
 = \K- \text{span of} ~ \left\{ \zeta(k_1, a_1/q_1) \cdots \zeta(k_r, a_r/q_r) ~:~ 1 \leq a_i \leq  q_i, 
 (a_i, q_i)=1, 1 \le i \le r \right\}.
 $$
 We observe that $V_{\K,\vec{k}}(\vec{q}) = \prod_{i=1}^r V_{\K,k_i}(q_i)$ and 
we define 
$V^+_{\K,\vec{k}}(\vec{q})= \prod_{i=1}^r V^+_{\K,k_i}(q_i)$. 

\end{defn}

In 2020, Hamahata proved the following theorem.

\begin{thm}{\rm (Hamahata \cite{YH})}
Let $q_1, \cdots , q_r $ be pairwise co-prime integers, $k_1,\cdots,k_r>1$ be positive integers
and $\vec{k}$, $\vec{q}$ be as in \defref{gCM}. If $\K$ is a number field such that
$\K \cap \Q(\zeta_{q_1\cdots q_{r}} )^+ = \Q$, then 
$$
\dim_{\K} V_{\K, \vec{k}}(\vec{q}) 
~\geq~
2^{-r}\prod_{i=1}^{r} \varphi(q_{i} ).
$$
\end{thm}

 Here we have the following extensions of Hamahata's theorem.
 
\begin{thm}\label{thm9}
Let $q_{t, j} \ge 1$ be integers for $1 \le j \le \ell$ and $1 \le t \le r$. Set $q_t = \prod_{j = 1}^ {\ell} q_{t, j}$
and $\vec{q_j} = ( q_{1,j}, \cdots ,  q_{r,j})$. 
Assume $q_1, \cdots , q_r $ are pairwise co-prime integers, $k_1,\cdots,k_r>1$ be positive integers
and $\vec{k}$, $\vec{q}$ be as in \defref{gCM}.
If $\K$ is a number field such that $\K \cap \Q(\zeta_{q_1 \cdots q_{r}})^+= \Q$, then 
$$
\dim_{\K} \left({\sum}_{j=1}^{\ell} V_{\K, \vec{k}}^+(\vec{q_j}) \right)
=
\begin{cases}
2^{-r}\sum_{j=1}^{\ell} \prod_{t=1}^{r} \varphi(q_{t , j}) & \text{  when at least one }k_t  \text{ is odd}, \\[2mm]
2^{-r}\sum_{j=1}^{\ell}\prod_{t=1}^{r} \varphi(q_{t, j}) - \ell+1 &    \text{  when all }k_t \text{ are even}.
 \end{cases}
$$ 
\end{thm}

\begin{thm}\label{thm10}
Let $r,\ell_1,\dots,\ell_r$ and $q_{t, j} \ge1$ be positive integers for $1 \le j \le \ell_t$ and $1 \le t \le r$. Set $q_t = \prod_{j = 1}^{\ell_t} q_{t, j}$ and $\vec{q_j} = ( q_{1,j}, \cdots , q_{r,j})$. Assume $q_1, \cdots , q_r $ are pairwise co-prime integers, $k_1,\cdots,k_r>1$ be positive integers
and $\vec{k}$ be as in \defref{gCM}.
If $\K$ is a number field such that $\K \cap \Q(\zeta_{q_1 \cdots q_{r}})^+= \Q$, then 
$$
\dim_{\K} \left(\prod_{t=1}^r{\sum}_{j=1}^{\ell_t} V_{\K, k_t}^+(q_{t, j}) \right)
= 2^{-r}\prod_{t=1}^{r} \left(\sum_{j=1}^{\ell_t}\varphi(q_{t, j}) - (\ell_t - 1)(1+(-1)^{k_t})\right)
$$ 
\end{thm}

\smallskip

\begin{rmk}
Let $[q_1, \cdots, q_{\ell}]$ denote the least common multiple of $q_1,\cdots, q_{\ell}$.
The number fields
$\K(\zeta_{\varphi(q_1)\cdots \varphi(q_{\ell})})$ and $\Q(\zeta_{q_1\cdots q_{\ell} })^{+}$
can be replaced by $\K(\zeta_{[\varphi(q_1),\cdots, \varphi(q_{\ell})]} )$ and
$\Q(\zeta_{[q_1,\cdots, q_{\ell} ]})^{+}$ respectively in 
Theorems \ref{thm1}, \ref{thm2}, \ref{cor1}, \ref{thm9} and Corollary \ref{cor3}.
\end{rmk}

\smallskip
 The article is organized as follows: in \S2 we list required results
needed for our proofs,  
in \S3 we derive main propositions and a corollary 
which are extensions of Okada's theorem 
(as well as extensions of Murty-Saradha) and a theorem of Hamahata.
Finally in the last section, we complete the proofs
 of Theorems \ref{thm1}, \ref{thm2},  \ref{cor1}, \ref{thm9}, \ref{thm10} and Corollary \ref{cor3}.

\medskip 
 
\section{\large{Preliminaries}}

\smallskip 

In this section, we fix some notations and state the results which 
will be used in the proofs of the main theorems.  When $p$ is an odd prime 
number, it is a result of Chowla \cite{chowla} that the set of numbers 
$$
\{ \cot(2\pi a/p) ~~|~~ 1 \leq a \leq (p-1)/2 \}
$$ 
are linearly independent over $\Q$.
This result was reproved by various authors (see for instance \cite{ HH, JH}).
In 1981, Okada \cite{TO} (see also Wang \cite{KW}) extended Chowla's theorem to 
natural numbers $q>2$. In the same article,
he also considered higher order derivatives of cotangent function.
More precisely, Okada~\cite{TO} proved the following theorem.
In order to state the theorem, let us denote $\frac{d^{k-1}}{dz^{k-1}} ( \cot z)|_{z = z_0}$ 
by $\cot^{(k-1)}(z_0)$.

\begin{thm}\label{okada} 
 Let $k$ and $q$ be positive integers with $k > 0$ and $q > 2$. Let $T$ be a
set of $\varphi(q)/2$ representatives modulo $q$ such that  $T \cup ( -T )$ is a complete set
of co-prime residues modulo $q$. Then the set of real numbers
$\{ \cot^{(k-1)}(\pi a/q) ~~|~ ~ a \in T \}$
is linearly independent over $\Q$.
\end{thm}
Using Galois theory, Girstmair \cite{ KG} gave an alternate proof 
for $\Q$ linear independence of derivatives of cotangent function. 
Murty-Saradha \cite{MS1} noticed that Okada's result can be extended 
to any number field $\K$ provided $\K \cap \Q(\zeta_q) =\Q$.
We note that the condition $\K \cap \Q(\zeta_q) =\Q$
in Murty-Saradha's result can be replaced by $\K \cap \Q(\zeta_q)^+ =\Q$.
Recently Hamahata \cite{YH} derived the following multi-dimensional
generalization of the above result.
 
\begin{thm}\label{hamahata}
Let $k_1, \cdots, k_r \geq 1$ be natural numbers and $q_1, \cdots, q_r >2 $ be 
co-prime natural numbers. For $1 \le i \le r$, let $T_i$ be a set of $\varphi(q_i)/2$ 
representatives modulo $q_i$ such that $T_i \cup (-T_i)$ 
is a complete set of co-prime residues modulo $q_i$. 
Set $q= q_1 \cdots q_r$. If $\K$ is a number field with $\K \cap \Q(\zeta_q) =\Q$, then
the $\varphi(q)/2^r$ numbers
\begin{equation*}
\prod_{i=1}^r \cot^{(k_{i}-1)}(\pi a_{i}/q_{i}),\quad a_{i}\in T_{i},\quad i =1,\cdots,r,
\end{equation*}
are linearly independent over $\K$.
\end{thm}
In the next section, we extend Theorems \ref{okada}, \ref{hamahata} for a finite set 
of pairwise co-prime natural numbers
$q_1, \cdots, q_{\ell} >2$. For this, we need to introduce linearly disjoint number fields.

\begin{defn} 
Let $\K$ and $\F$ be algebraic extensions of a field $\L$. The fields $\K, \F$ 
are said to be linearly disjoint over $\L$ if every finite subset 
of $\K$ that is $\L$ linearly independent is also $\F$  linearly independent. 
\end{defn}
We have the following equivalent criterion for linearly disjoint fields.

\begin{thm}\label{kf} \cite[Ch 5, Thm 5.5]{MC} 
Let $\K$ and $\F$ be algebraic extensions of a field $\L$. Also let
at least one of $\K, \F$ is separable and one (possibly the same) is normal. 
Then $\K$ and $\F$ are linearly disjoint over $\L$ if and only if $\K \cap \F = \L$.
\end{thm}

\smallskip

\smallskip

\section{\large{Requisite Propositions}}

We first relate the vector spaces $V_{\Q,k}^+(q)$ to $\Q$-vector subspaces of cyclotomic fields, 
via the cotangent values mentioned in the previous section.

Indeed, since $\pi\cot(\pi z)$ is the logarithmic derivative of 
$\sin(\pi z) = z\prod_{m=1}^\infty(1-z^2/m^2)$, one checks for $k>1$ that
$$\pi^{k}\cot^{(k-1)}(\pi z) = (-1)^{k-1}(k-1)!\left(\frac{1}{z^{k}} 
+ \sum_{m=1}^\infty\left(\frac{1}{(z+m)^{k}}+\frac{1}{(z-m)^{k}}\right)\right)
$$
and, evaluating at $z=a/q$ with $a/q\notin\Z\cup\{\infty\}$, it follows
\begin{equation}\label{formulecotg}
\zeta (k,a/q)+(-1)^{k}\zeta (k,1-a/q) = \frac{\pi^k(-1)^{k-1}}{(k-1)!}\cot^{(k-1)}(\pi a/q).
\end{equation}
Note that for $z \notin \Z$, we have (see \cite{NR})
\begin{equation}\label{dcot}
\cot^{(k-1)} (\pi z) 
= \sum_{a, b \ge 0 \atop a + 2b = k} 
\beta_{a, b}^{(k)} \cot^a \pi z 
~(1 + \cot^2 \pi z)^b,
\end{equation}
where $\beta_{a, b}^{(k)} \in \Z$ and also
\begin{equation}\label{ecot}
-i \cot \frac{\pi a}{q} 
~=~ 
\frac{\zeta^{a}_{q} + 1}{\zeta^{a}_{q} - 1} 
~\in~
\Q(\zeta_{q}),
\end{equation} 
where $ i = \sqrt{-1}$.  We then observe that
\begin{eqnarray*}
(i\pi)^{-k} V_{\Q,k}^+(q) 
&=& 
\sum_{\stackrel{1\le a\le q/2}{(\scriptscriptstyle a,q)=1}} 
\Q ( i^k\cot^{(k-1)}(\pi a/q)) 
~ \subset \Q(\zeta_q) \cap (i^k\R) \\
\phantom{m}\text{and}\phantom{m}
&&
(\zeta_q - \zeta_q^{-1})^k i^k\cot^{(k-1)}(\pi a/q)  \in 
\Q(\zeta_q)^{+}.
\end{eqnarray*}
Since the dimension of the $\Q$-vector space $\Q(\zeta_q)\cap(i^k\R)$ is 
$\varphi(q)/2$, and $\dim_{\Q}V^+_{\Q,k}(q) = \varphi(q)/2$ 
by Okada's theorem \ref{okada}, we have 
\begin{equation}\label{pat}
\Q(\zeta_q)\cap(i^k\R) = (i\pi)^{-k}V_{\Q,k}^+(q)
\end{equation}
for $q> 2$. It is easy to see that the above relation also holds for $q=1, 2$.

\smallskip

\begin{lem}\label{leminter}
If $k>1$, $q$, $q_1$, $q_2\ge1$ are positive integers and $d=(q_1, q_2)$, then 
\begin{equation*}
(i\pi)^{-k}V^+_{\Q,k}(q) = \Q(\zeta_q)\cap(i^k\R), ~
V^+_{\Q,k}(q)\subset V^+_{\Q,k}(q_1) \text{ if } q|q_1
\text{  and} \phantom{m}
V_{\Q, k}^+(q_1) \cap V_{\Q, k}^+(q_2) = V_{\Q,k}^+(d).
\end{equation*}
\end{lem}
\begin{proof}
Identity \eqref{pat} asserts precisely the first statement
\begin{equation*}
\Q(\zeta_q)\cap(i^k\R) = (i\pi)^{-k}V_{\Q,k}^+(q). 
\end{equation*}
If $q|q_1$, we have $\Q(\zeta_q)\subset \Q(\zeta_{q_1})$. We also 
know from Galois theory that
\begin{equation}\label{gal}
\Q(\zeta_{q_1}) \cap \Q(\zeta_{q_2}) = \Q(\zeta_{d}).
\end{equation}
Then the second and the third results follow by intersecting both sides of 
the inclusion and Equation~\eqref{gal} with $i^k\R$ and applying \eqref{pat}.
\end{proof}

\begin{prop}\label{propinter}
For an integer $k>1$ and pairwise co-prime positive integers $q_1,\cdots ,q_\ell$, 
the kernel of the surjective map
$$\begin{matrix}
\oplus_{j=1}^\ell V_{\Q, k}^+(q_j) &\longrightarrow &\sum_{j=1}^{\ell}V_{\Q, k}^+(q_j)\\
(x_1,\cdots,x_\ell) &\longmapsto &x_1+\cdots+x_\ell
\end{matrix}$$
is $0$ if $k$ is odd. When $k$ is even, the kernel of the above map is the
$\Q$-vector space $V_{\Q, k}^+(1)^{\ell-1}$ of dimension $\ell-1$ which is parametrised as
$$
(z_1,\cdots,z_{\ell-1})\in\Q^{\ell-1} \longmapsto(z_1\zeta(k,1),\cdots,z_{\ell-1}\zeta(k,1),
-(z_1+\cdots+z_{\ell-1})\zeta(k,1)) \in \oplus_{j=1}^\ell V_{\Q, k}^+(q_j).
$$
\end{prop}
\begin{proof}
The kernel of the map consists of elements $(x_1, \cdots, x_{\ell})$ where
$x_j\in V_{\Q, k}^+(q_j)\cap \sum_{t\not=j}V_{\Q, k}^+(q_t)$
for all $j=1,\cdots,\ell$. 
By Lemma \ref{leminter}, $\sum_{t\not=j}V_{\Q, k}^+(q_t)\subset V_{\Q, k}^+(\prod_{t\not=j}q_t)$ 
and, since $q_j$ is co-prime to $\prod_{t\not=j}q_t $, 
it also implies $x_j\in V_{\Q, k}^+(1)$ for $j=1,\cdots,\ell$.
\end{proof}

We will need a multivariate extension of this result. To do this, we first prove another lemma.

\begin{lem}\label{leminterm}
Let $k_1,\cdots,k_r>1$ be integers and 
$q_{1,1},\cdots ,q_{r,1}, q_{1,2},\cdots ,q_{r,2}$ be positive integers which are 
pairwise co-prime. Set $\vec{k}=(k_1,\cdots,k_r)$ and 
$\vec{q}_j=(q_{1,j},\cdots,q_{r,j})$ for $j=1,2$, then
$$
V_{\Q, \vec{k}}^+(\vec{q}_{1}) \cap V_{\Q, \vec{k}}^+(\vec{q}_{2}) 
=
V_{\Q,\vec{k}}^+(\vec{1}).
$$ 
\end{lem}
\begin{proof}
Applying \lemref{leminter}, it is easy to see that 
$$
V_{\Q,\vec{k}}^+(\vec{1}) 
\subseteq 
V_{\Q, \vec{k}}^+(\vec{q}_{1}) \cap V_{\Q, \vec{k}}^+(\vec{q}_{2}).
$$
Set $q_j=\prod_{t=1}^rq_{t, j}$ and $k=\sum_{t=1}^rk_t$. 
Again Lemma \ref{leminter} shows that 
$(i\pi)^{-k}V_{\Q, \vec{k}}^+(\vec{q}_{j})\subset \Q(\zeta_{q_j})$ for $j=1,2$ and, 
since $q_1$ and $q_2$ are relatively prime, the intersection is contained in $\Q$. 
Now for $j=1$ or $2$, if some rational number $\alpha$ belongs to 
$(i\pi)^{-k}V_{\Q, \vec{k}}^+(\vec{q}_{j})$, we have for each $t =1, \cdots, r$,  
a linear relation over $\Q(\zeta_{\prod_{h\not=t}q_{h,j}})$ between $\alpha$ 
and the elements of a basis (over $\Q$) of $(i\pi)^{-k_t}V_{\Q,k_t}^+(q_{t, j})$. 
The field $\Q(\zeta_{\prod_{h\not=t}q_{h,j}})$ is linearly disjoint from 
$\Q(\zeta_{q_{t, j}})$ because $q_{t, j}$ is co-prime to $\prod_{h\not=t}q_{h,j}$. 
Thus the above linear relation over $\Q(\zeta_{\prod_{h\not=t}q_{h,j}})$ 
implies a linear relation over $\Q$, which involves $\alpha$, since the 
elements of a basis of $(i\pi)^{-k_t}V_{\Q,k_t}^+(q_{t, j})$ are linearly 
independent over $\Q $. Hence $\alpha\in(i\pi)^{-k_t}V_{\Q,k_t}^+(q_{t, j})$
for all $t=1,\cdots, r$ and $j=1, 2$.  Applying \lemref{leminter}, we see that 
$\alpha \in (i\pi)^{-k_t}V_{\Q, k_t}^+(1)$. Note that
$(i\pi)^{-k_t}V_{\Q, k_t}^+(1)$ is $0$ if $k_t$ is odd and 
$\Q$ if $k_t$ is even. Therefore, if $\alpha\ne0$ then all $k_t$ must be even and 
$$
V_{\Q, \vec{k}}^+(\vec{q}_{1}) \cap V_{\Q, \vec{k}}^+(\vec{q}_{2}) 
\subseteq \prod_{t=1}^r V^+_{\Q,k_t}(1) =
\begin{cases}
0 & \text{  when at least one }k_t  \text{ is odd}, \\[2mm]
\Q\prod_{t=1}^r \zeta(k_t, 1) &    \text{  when all }k_t \text{ are even},
\end{cases}
$$
which is equal to $V_{\Q,\vec{k}}^+(\vec{1})$.
\end{proof}

\begin{prop}\label{propinterm}
For integers $k_1,\cdots,k_r>1$ and pairwise co-prime positive integers 
$q_{t,j}$, $t=1,\cdots,r$, $j=1,\cdots,\ell$, we define 
$\vec{k}=(k_1,\cdots,k_r)$ and $\vec{q}_j=(q_{1,j},\cdots,q_{r,j})$. 
Then the kernel of the surjective map
$$\begin{matrix}
\oplus_{j=1}^\ell V_{\Q, \vec{k}}^+(\vec{q}_{j}) &\longrightarrow &\sum_{j=1}^{\ell}V_{\Q, \vec{k}}^+(\vec{q}_{j})\\
(x_1,\cdots,x_\ell) &\longmapsto &x_1+\cdots+x_\ell
\end{matrix}$$
is $0$ if at least one $k_t$ is odd and is the $\Q$-vector space 
$\left(V_{\Q, \vec{k}}^+(\vec{1})\right)^{\ell-1}$ of dimension $\ell -1$ parametrised as
$$\begin{matrix}
\hfill\Q^{\ell-1} &\longrightarrow &\oplus_{j=1}^\ell V_{\Q, \vec{k}}^+(\vec{q}_{j})\hfill\\
(z_1,\cdots,z_{\ell-1}) &\longmapsto &\big(z_1\prod_{t=1}^r\zeta(k_t,1),\cdots,z_{\ell-1}\prod_{t=1}^r\zeta(k_t,1),-(z_1+\cdots+z_{\ell-1})\prod_{t=1}^r\zeta(k_t,1)\big)
\end{matrix}$$
if all $k_1,\cdots,k_r$ are even.
\end{prop}
\begin{proof}
The kernel of the map consists of elements $(x_1, \cdots, x_{\ell})$, where
$$
x_j\in V_{\Q, \vec{k}}^+(\vec{q}_{j})\cap \sum_{h\not=j} V_{\Q, \vec{k}}^+(\vec{q}_{h})
$$
for $j=1,\cdots,\ell$.  By Lemma \ref{leminter}, we have 
$\sum_{h\not=j}V_{\Q, \vec{k}}^+(\vec{q}_{h})\subset V_{\Q, \vec{k}}^+(\vec{Q})$,
 where $\vec{Q}=(Q_1,\cdots,Q_r)$ with $Q_t=\prod_{h\not=j}q_{t,h}$. 
 Since $q_{1,j},\cdots,q_{r,j},Q_1,\cdots,Q_r$ are pairwise co-prime, 
 Lemma \ref{leminterm} implies that $x_j=0$ if at least one 
 $k_t$ is odd and $x_j\in \prod_{t=1}^r V_{\Q, k_t}^+(1) = \Q\prod_{t=1}^r\zeta(k_t,1)$ 
 if all $k_1,\cdots,k_r$ are even, for $j=1,\cdots,\ell$.
\end{proof}

\smallskip

For $1 \leq j \leq \ell $ and $q_j >2$, consider the sets
\begin{gather}\label{def-t}
S_{q_j} = \{ 1 < a_j < q_j /2 ~ : ~ (a_j,q_j) = 1 \}
\phantom{m}\text{and}\phantom{m}
T_{q_j} = \{ 1\leq a_j < q_j /2 ~ : ~ (a_j,q_j) = 1 \},\nonumber\\
U_{j}(k)
~=~
\begin{cases}
T_{q_{j}} & \text{  when }k \text{ is odd}, \\
T_{q_{1}}&    \text{  when }j=1 \text{ and }k \text{ is even},\\
S_{q_{j}}&    \text{  when }j \ne 1 \text{ and }k \text{ is even}.
\end{cases}
\end{gather}
We have the following proposition about linear independence of
derivatives of cotangent function evaluated at certain rational points. 
It forms a basis for the space $\sum_{j=1}^\ell V^+_{\K,k}(q_j)$.

\begin{prop}\label{prop1}
For $1 \le j \le \ell$, let $q_j >2$  be pairwise co-prime natural numbers and 
$k > 1$ be an integer. If $\K$ is a number field such that 
$\K \cap \Q(\zeta_{q_1\cdots q_{\ell} })^+ = \Q$, then the set of real numbers
$$
\bigcup_{ 1 \leq j \leq \ell} \left\{ \cot^{(k-1)} \left(\frac{\pi a_j}{q_j}\right) ~ : ~ a_j \in U_j(k)\right\}
$$
is linearly independent over $\K$.
\end{prop}

\begin{proof}
By \eqref{formulecotg}, the numbers 
$\bigcup_{ 1 \leq j \leq \ell} \left\{ \pi^k\cot^{(k-1)} \left(\frac{\pi a_j}{q_j}\right) ~ : ~ a_j \in T_{q_j}\right\}$ 
span $\sum_{\j=1}^\ell V_{\K,k}^+(q_j)$. 
When $k$ is even and $q > 2$ is any integer, using \eqref{eq2} and replacing $k$ by $2k$,
we have
$$
\zeta(2k,1) 
~=~
 \alpha\sum_{1 \leq a < q/2 \atop (a,q)=1 }\left( \zeta(2k, a/q) + \zeta(2k, 1-a/q) \right)
~=~ 
\alpha \sum_{1 \leq a < q/2 \atop (a,q)=1 }\pi^{2k}\cot^{(2k-1)} \left(\frac{\pi a}{q}\right),
$$
where $\alpha$  is some non zero rational number. It implies that, removing 
one arbitrary element from any $\ell -1$ sets of the form
$\left\{ \pi^{2k}\cot^{(2k-1)} \left(\frac{\pi a_j}{q_j}\right) ~ : ~ a_j \in T_{q_j}\right\}$
in the union
$$
\bigcup_{ 1 \leq j \leq \ell} \left\{ \pi^{2k}\cot^{(2k-1)} \left(\frac{\pi a_j}{q_j}\right) ~ : ~ a_j \in T_{q_j}\right\}
$$ 
does not change the space they span, 
that is $\sum_{\j=1}^\ell V_{\K,2k}^+(q_j)$. Whatever the value of $k>1$ is, 
the number of generators left is equal to the dimension of this span. 
Thus these generators must be linearly independent.
\end{proof}

\smallskip

We have the following corollary of \propref{prop1} which gives a basis
of the $\K$ vector space $\sum_{j=1}^\ell V^+_{\K,\vec{k}}(\vec{q_j})$. 

\begin{cor}\label{coha}
	Let $q_{t,j} >  2, T_{q_{t,j}} $ be defined as before \propref{prop1} 
	and $q_t= \prod_{j = 1}^ {\ell} q_{t,j}$. 
	Assume $q_1, \cdots , q_r $ are pairwise co-prime and $k_1,\cdots,k_r$ 
	be positive integers. Suppose that
	$$U_{j}(\vec{k}) =
	\begin{cases}
		\prod_{t=1}^rT_{q_{t,j}} & \text{  when one of the }k_h \text{ is odd}, \\
		\prod_{t=1}^rT_{q_{t,1}} &   \text{  when }j=1 \text{ and all the }k_h\text{ are even},\\
		\prod_{t=1}^rT_{q_{t,j}}\setminus \{(1, \cdots, 1)\}	&    \text{  when }j \ne 1 \text{ and all the }k_h \text{ are even}.
	\end{cases}$$
	If $\K$ is a number field such
	that $\K \cap \Q(\zeta_{q_1 \cdots q_{r}})^+= \Q$,
	then the set of real numbers
	\begin{equation*}
		\bigcup_{j=1}^{\ell} \left\{ ~\prod_{t=1}^r \cot^{(k_t - 1)}\left(\frac{\pi a_{t}}{q_{t,j}}\right) 
		~:~ (a_1, \cdots,  a_r) \in U_{j}(\vec{k})  \right\}
	\end{equation*}
	is linearly independent over $\K$.
\end{cor}

\begin{proof}
	We will use induction on $r$ to complete the proof. When $r=1$ and $k_1=1$, the 
	result follows from the work of the first and second author \cite{GK}. 
	If $r=1$ and $k_1 >1$, result follows from \propref{prop1}. 
	Now assume that the corollary is true for any natural number strictly less than
	$r$. Set $q = q_1 \cdots q_r$.
	Suppose that there exist rational numbers $\alpha_{j,a_{1}, \cdots, a_{r}}$
	such that
	$$
	\sum_{j=1}^{\ell} \sum_{a_{1} }\cdots \sum_{a_{r}}
	\alpha_{j,a_{1}, \cdots, a_{r}} 
	\prod_{t=1}^r  i^{k_t} \cot^{(k_t - 1)}\left(\frac{\pi a_{t}}{q_{t,j}}\right) = 0,
	$$
	where $(a_1, \cdots,  a_r)$ runs over the elements of $ U_{j}(\vec{k})$.
	This implies that
	\begin{equation*}
		\sum_{j=1}^{\ell} \sum_{a_{r} } 
		\left(\sum_{a_{1}} \cdots \sum_{a_{r-1} 
		} \alpha_{j,a_{1}, \cdots, a_{r}} 
		\prod_{t=1}^{r-1}  i^{k_t} \cot^{(k_t - 1)}
		\left(\frac{\pi a_t}{q_{t,j}}\right)
		\right)
		i^{k_r}
		\cot^{(k_{r}-1)}
		\left(\frac{\pi a_{r}}{q_{r,j}}\right) = 0.
	\end{equation*}
	Since
	$\Q(\zeta_{q_1 \cdots q_{r-1} }) \cap \Q(\zeta_{q_r} ) = \Q$,
	\propref{prop1} implies
	$$
	\sum_{a_{1} } \cdots \sum_{a_{r-1}} 
	\alpha_{j,a_{1}, \cdots, a_{r}}
	\prod_{t=1}^{r-1}  i^{k_t}  \cot^{(k_t - 1)}\left(\frac{\pi a_{t}}{q_{t,j}}\right) = 0,
	$$
	for $ 1 \le j \le \ell $ and $(a_1, \cdots,  a_r) \in U_{j}(\vec{k})$. By induction hypothesis,
	$\alpha_{j,a_{1}, \cdots, a_{r}} = 0 $ for all $1 \le j \le \ell$ and $(a_1, \cdots,  a_r) \in U_{j}(\vec{k})$.
	Therefore,
	the set of real numbers
	\begin{equation*}
		\bigcup_{j=1}^{\ell} \left\{ ~\prod_{t=1}^r i^{k_t}  \cot^{(k_t - 1)}\left(\frac{\pi a_{t}}{q_{t,j}}\right) 
		~:~ (a_1, \cdots,  a_r) \in U_{j}(\vec{k}) \right\}
	\end{equation*}
	is linearly independent over $\Q$.
	It implies that
	\begin{equation*}
		\bigcup_{j=1}^{\ell} \left\{ ~\prod_{t=1}^r i^{k_t} 
		(\zeta_q - \zeta^{-1}_{q})^{k_t}  \cot^{(k_{t}-1)}\left(\frac{\pi a_{t}}{q_{t,j}}\right) 
		~:~ (a_1, \cdots,  a_r) \in U_{j}(\vec{k}) \right\}
	\end{equation*}
	is linearly independent over $\Q$. As in \propref{prop1}, the numbers inside the products belong to $\Q(\zeta_q)^+$
	and by given hypothesis, $\K$ and $\Q(\zeta_q)^+$ are linearly disjoint.
	Hence the numbers in the union are linearly independent over $\K$ as well. It implies that
	\begin{equation*}
		\bigcup_{j=1}^{\ell} \left\{ ~\prod_{t=1}^r  \cot^{(k_t -1)}\left(\frac{\pi a_{t}}{q_{t,j}}\right) 
		~:~ (a_1, \cdots,  a_r) \in U_{j}(\vec{k}) \right\}
	\end{equation*}
	is linearly independent over $\K$.
\end{proof}

\begin{rmk}
Let $\chi$ be a character modulo $q$. If $\chi$ has the same parity as $k$ (\emph{i.e.} $\chi(-1)=(-1)^k$), then
\begin{equation*}
\frac{\pi^k(-1)^{k-1}}{q^k(k-1)!}\sum_{1 \le a < q/2, \atop{(a,q)=1 }}\chi(a)\cot^{(k-1)}(\pi a/q) 
= q^{-k}\sum_{(a,q)=1}\chi(a)\zeta (k, a/q) = L(k,\chi).
\end{equation*}
If $\chi$ and $k$ have different parity (\emph{i.e.} $\chi(-1)=(-1)^{k-1}$),
then 
$$
\sum_{1 \le a < q, \atop{ (a,q)=1 }}\chi(a)\cot^{(k-1)}(\pi a/q) = 0
$$ 
since $\cot^{(k-1)}(-\pi a/q) = (-1)^k\cot^{(k-1)}(\pi a/q)$.
\end{rmk}

\smallskip

\section{\large{Proofs of the Main Theorems}}

\subsection{Proof of Theorem \ref{cor1}}
 Since, by hypothesis, the fields $\K$ 
and $\Q(\zeta_{q_1\cdots q_\ell})^+$ are linearly disjoint, it suffices 
to determine the dimension of the $\Q$-vector space 
$\sum_{j=1}^\ell (i\pi)^{-k}V_{\Q,k}^+(q_j)$. But, this is the dimension 
of $\oplus_{j=1}^\ell (i\pi)^{-k}V_{\Q,k}^+(q_j)$ minus the dimension 
of the kernel of the map shown in Proposition \ref{propinter}. 
This gives the first equality, the remaining ones coming from 
$$
\dim_{\Q}(V_{\Q,k}^+(q_j)) = \varphi(q_j)/2
\phantom{m}\text{and}\phantom{m} 
\dim_{\Q}(V_{\Q,k}^+(1)) = \dim_{\Q}(V_{\Q,k}^+(2)) = \frac12(1+(-1)^k).
$$
\qed

\subsection{Proof of Corollary \ref{cor3}}
The lower bounds follow directly from Theorem \ref{cor1}. Recall that for $k>1$ and $q >2$, we have
\begin{equation}\label{eqdsprcor}
\zeta(k) \prod_{p | q \atop{p \text{ prime} } } (1- p^{-k})  
= \frac{1}{ q^{k}}\sum_{1 \leq a< q \atop (a, q)=1 } \zeta(k, a/ q)
\end{equation}
and
$$V_{\K, k}(q) = \K- \text{span of} ~ \{ \zeta(k,a/q) : 1 \leq a < q, (a,q)=1\}.
$$
The number of generators of $\sum_{j=1}^{ \ell} V_{\K, k}(q_j)$ is $\sum_{j=1}^{\ell} \varphi(q_j)$. But, by \eqref{eqdsprcor}, there are at least $\ell-1$ independent linear relations between them.
Hence
$$
\dim_{\K} \left(\sum_{j=1}^{ \ell} V_{\K, k}(q_j) \right) 
\leq 
~\sum_{j=1}^{\ell} \varphi(q_j) - (\ell - 1). 
$$
\qed

\subsection{Proof of \thmref{thm1}}
Let $\F = \K (\zeta_{\varphi(q_1)\cdots\varphi(q_\ell)})$.
By \eqref{dry}, the dimension of the $\F$ vector space generated by the elements of 
$X_{q_1,2k+1,o}\cup\cdots\cup X_{q_\ell,2k+1,o}$ is equal to the dimension of 
$\sum_{j=1}^\ell V_{\F ,2k+1}^+(q_j)$. Since, by hypothesis, 
$\F \cap\Q(\zeta_{q_1\cdots q_\ell})^+=\Q$, 
Theorem \ref{cor1} shows that this dimension is equal to the number of 
elements of $X_{q_1,2k+1,o}\cup\cdots\cup X_{q_\ell,2k+1,o}$, which 
must therefore be linearly independent over $\F$.

\smallskip

Alternatively we can argue as follows.
For $1 \leq j \leq \ell$, let $D_j = \{\chi_j \mod{q_j} ~ | ~ \chi_j(-1) = -1\}$
be the set of odd characters modulo $q_j$. 
Suppose that there exist $\alpha_{\chi_{j}} \in \K(\zeta_{\varphi(q_1)\cdots\varphi(q_{\ell})})$ 
for $\chi_j \in D_j, 1 \leq j \leq \ell$ such that 
\begin{equation}\label{der}
\sum_{ 1 \leq j \leq \ell} \sum_{\chi_j \in D_j} \alpha_{\chi_{j}} L(2k+1,\chi_j) = 0.
\end{equation}
Substituting (see \cite{NR})
\begin{equation*}
L(2k+1, \chi_j)  = \frac{\pi^{2k+1}}{ (2k)! ~ ~ q^{2k+1}_j} 
\sum_{a_j \in T_{q_j}} \chi_j(a_j) \cot^{(2k)} (\frac{\pi a_j}{q_j}), 
\end{equation*}
for $\chi_j \in D_j,1 \leq j \leq \ell$ in \eqref{der}, we obtain
\begin{align}
\begin{split}\label{lkchi}
\sum_{ 1 \leq j \leq \ell} \sum_{a_j \in T_{q_j}} \frac{1}{q^{2k+1}_j} \left(\sum_{\chi_j \in D_j} 
\alpha_{\chi_{j}}\chi_j(a_j)\right) \cot^{(2k)} (\frac{\pi a_j}{q_j}) &= 0.
\end{split}
\end{align}
Here $T_{q_j}$ is as defined in \eqref{def-t}. By given hypothesis, we have
$\K(\zeta_{\varphi(q_1)\cdots\varphi(q_{\ell})}) \cap \Q(\zeta_{q_1 \cdots q_{\ell}})^+ = \Q$.
 It then follows from \propref{kf} that $\K(\zeta_{\varphi(q_1)\cdots\varphi(q_{\ell})})$ 
and $\Q(\zeta_{q_1 \cdots q_{\ell}})^+$ are linearly disjoint over $\Q$. 
Note that the coefficients of $\cot^{(2k)} (\frac{\pi a_j}{q_j})$'s in \eqref{lkchi}
belong to $\K(\zeta_{\varphi(q_1)\cdots\varphi(q_{\ell})})$ and hence \propref{prop1}
implies that
$$
\sum_{\chi_j \in D_j} \alpha_{\chi_{j}} \chi_j(a_j) ~=~~ 0
$$
for $a_j \in T_{q_j} ,1 \leq j \leq \ell$. Since all the characters in the set $D_j, 1 \leq j \leq \ell$ are 
of same parity, it follows that 
$$
\sum_{\chi_j \in D_j} \alpha_{\chi_{j}} \chi_j(a_j) = 0
$$ 
for $a_j \in (\Z/q_j\Z)^{\times}  , 1 \leq j \leq \ell$. It then follows
from linear independence of characters that
$\alpha_{\chi_{j}} = 0$ for $\chi_j \in D_j, 1 \leq j \leq \ell$.
This completes the proof of \thmref{thm1}.
\qed

\smallskip

\subsection{Proof of \thmref{thm2}}
Let $\F = \K (\zeta_{\varphi(q_1)\cdots\varphi(q_\ell)})$.
By \eqref{dry}, the dimension of the $\F$ vector space generated by the elements of 
$\{\zeta(2k)\}\cup X_{q_1,2k,e}\cup\cdots\cup X_{q_\ell,2k,e}$ is equal 
to the dimension of $\sum_{j=1}^\ell V_{\F, 2k}^+(q_j)$. Since, by hypothesis, 
$\F \cap \Q(\zeta_{q_1\cdots q_\ell})^+=\Q$, it follows from
Theorem \ref{cor1} that the dimension of $\sum_{j=1}^\ell V_{\F, 2k}^+(q_j)$ is equal to the number 
of elements of $\{\zeta(2k)\}\cup X_{q_1,2k,e}\cup\cdots\cup X_{q_\ell,2k,e}$, which must therefore be linearly independent over $\F$.

\smallskip

Alternatively we can argue as follows.
 Let $C_1 = \{\chi_1 \mod{q_1} ~ | ~ \chi_1(-1) = 1\}$ and for $1 < j \leq \ell$,
 let $C_j = \{\chi_j \mod{q_j} ~ | ~ \chi_j(-1) = 1, \chi_j \ne 1\}$
be the set of non-trivial even characters modulo $q_j$. 
Suppose that there exist $\alpha_{\chi_{j}} \in 
\K(\zeta_{\varphi(q_1)\cdots\varphi(q_{\ell})})  $ for $\chi_j \in C_j, 1 \leq j \leq \ell$ such that 
\begin{equation}\label{der1}
	\sum_{ 1 \leq j \leq \ell} \sum_{\chi_j \in C_j} \alpha_{\chi_{j}} L(2k,\chi_j) = 0.
\end{equation}
Substituting (see \cite{NR})
\begin{equation}\label{odd-chi1}
	L(2k,\chi_j)  = \frac{-\pi^{2k}}{ (2k-1)! ~ ~ q^{2k}_j} \sum_{a_j \in T_{q_j}} \chi_j(a_j) \cot^{(2k-1)} (\frac{\pi a_j}{q_j}), 
\end{equation}
for $\chi_j \in C_j,1 \leq j \leq \ell$ in \eqref{der1}, we obtain
\begin{align}
	\begin{split}\label{lkchi1}
		\sum_{ 1 \leq j \leq \ell} \sum_{a_j \in T_{q_j}}\frac{1}{q^{2k}_j} 
		 \left(\sum_{\chi_j \in C_j} 
		\alpha_{\chi_{j}} \chi_j(a_j)\right) \cot^{(2k-1)} (\frac{\pi a_j}{q_j}) &= 0,
	\end{split}
\end{align}
where $T_{q_j}$ is as in \eqref{def-t}. For $a_j \in T_{q_j}, 1 \leq j \leq \ell$, let us denote
$$
 A_j(a_j)= \frac{1}{q^{2k}_j} \sum_{\chi_j \in C_j} 
\alpha_{\chi_{j}} \chi_j(a_j).
$$
From \eqref{lkchi1}, we have 
\begin{equation}\label{ngt}
\sum_{ 1 \leq j \leq \ell} \sum_{a_j \in T_{q_j}} A_j(a_j) \cot^{(2k-1)} (\frac{\pi a_j}{q_j}) 
~=~ 0.
\end{equation}
This implies that
\begin{eqnarray*}
\sum_{a_1 \in T_{q_1}} A_1(a_1) \cot^{(2k-1)} (\frac{\pi a_1}{q_1}) 
&+& 
\sum_{1 < j \leq \ell} A_j(1) \sum_{a_j \in T_{q_j}}  \cot^{(2k-1)} (\frac{\pi a_j}{q_j})  \\
&+& \sum_{ 1 < j \leq \ell} \sum_{a_j \in S_{q_j}} ( A_j(a_j) - A_j(1) ) \cot^{(2k-1)} (\frac{\pi a_j}{q_j}) 
=0,
\end{eqnarray*}
where $S_{q_j} = \{ 1 < a_j < q_j /2 ~ | ~ (a_j,q_j) = 1 \}$. As in the previous section, for $1 < j \leq \ell$, 
recalling 
$$
 \sum_{a_{j} \in T_{q_j}} \cot^{(2k-1)} (\frac{\pi a_j}{q_j})
 =
\frac{tr_j}{tr_1} \sum_{a_{1} \in T_{q_1}} \cot^{(2k-1)} (\frac{\pi a_1}{q_1}),
$$ 
we have
\begin{equation}\label{song}
\sum_{a_1 \in T_{q_1}} (A_1(a_1) + \sum_{ 1< j \leq \ell} \frac{ A_j(1)tr_j}{tr_1} )
\cot^{(2k-1)} (\frac{\pi a_1}{q_1}) 
+
 \sum_{ 1 < j \leq \ell} \sum_{a_j \in S_{q_j}} 
\left (A_j(a_j) - A_j(1) \right) \cot^{(2k-1)} (\frac{\pi a_j}{q_j}) = 0.
\end{equation}
As in \thmref{thm1}, given hypothesis implies that the fields $\K(\zeta_{\varphi(q_1)\cdots\varphi(q_{\ell})})$ 
and  $\Q(\zeta_{q_1 \cdots q_{\ell}})^+$ are linearly disjoint. 
Since the coefficients of $\cot^{(2k-1)} (\frac{\pi a_j}{q_j})$'s in \eqref{song}
belong to $\K(\zeta_{\varphi(q_1)\cdots\varphi(q_{\ell})})$, using \propref{prop1}, we obtain
\begin{equation}\label{tab}
	A_1(a_1) + \sum_{ 1< j \leq \ell} \frac{ A_j(1)tr_j}{tr_1}  = 0 \phantom{mm} \text{and} \phantom{mm}
	A_j(a_j) = A_j(1),
\end{equation}
where  $ a_1 \in T_{q_1}$ and $a_j \in S_{q_j} $ for $1 < j \leq \ell$. The second equality implies that
  \begin{align*}
 \frac{1}{q^{2k}_j}    \sum_{\chi_j \in C_j} 
  \alpha_{\chi_{j}}\chi_j(a_j) ~=~ A_j(1) \chi_0(a_j),
\end{align*}
where $\chi_0$ is the trivial character modulo $q_j$ and $a_j \in T_{q_j} $.
Since all the characters in the set $C_j, 1 \leq j \leq \ell$ are even, it follows that 
\begin{align*}
\frac{1}{q^{2k}_j} \sum_{\chi_j \in C_j} 
\alpha_{\chi_{j}} \chi_j(a_j) ~=~ A_j(1) \chi_0(a_j),
\end{align*}
where $a_j \in (\Z/q_j\Z)^{\times}, 1 < j \leq \ell$. Linear
independence of characters implies that  $A_j(1) = 0$  and
$\alpha_{\chi_{j}} = 0$ for $\chi_j \in C_j, 1 < j \leq \ell$.
Replacing $A_j(1) = 0$ in the first equality of \eqref{tab}, we get
$A_1(a_1)=0$ for $a_1\in T_1$. Arguing as above, we get
that $\alpha_{\chi_1} =0$ for $\chi_1 \in C_1$.
This completes the proof of \thmref{thm2}.

\smallskip

\subsection{Proof of \thmref{thm9}}
Set $k = \sum_{t=1}^r k_t$. 
As in the proof of \thmref{cor1}, by hypothesis, the fields $\K$ and 
$\Q(\zeta_{q_1\cdots q_r})^+$ are linearly disjoint, thus it suffices 
to determine the dimension of the $\Q$-vector space 
$\sum_{j=1}^\ell (i\pi)^{-k}V_{\Q,\vec{k}}^+(\vec{q}_j)$. But, this is 
the dimension of $\oplus_{j=1}^\ell (i\pi)^{-k}V_{\Q,\vec{k}}^+(\vec{q}_j)$ 
minus the dimension of the kernel of the map shown in 
Proposition \ref{propinterm}. Since 
$$
\dim_{\Q}\left( \oplus_{j=1}^\ell (i\pi)^{-k}V_{\Q,\vec{k}}^+(\vec{q}_j) \right) 
= 2^{-r}\sum_{j=1}^\ell\prod_{t=1}^r\varphi(q_{t,j}),
$$ 
the result follows from the same Proposition \ref{propinterm}.
\qed

\subsection{Proof of \thmref{thm10}}
The map
$$\begin{matrix}
\otimes_{t=1}^r\sum_{j=1}^{\ell_t}(i\pi)^{-k_t}V_{\Q,k_t}^+(q_{t, j}) 
&\longrightarrow &\prod_{t=1}^r\sum_{j=1}^{\ell_t}(i\pi)^{-k_t}V_{\Q,k_t}^+(q_{t,j})\\
\hfill x_1\otimes\cdots\otimes x_r &\longmapsto &x_1\cdots x_r\hfill
\end{matrix}
$$
is a bijection. Indeed, a nonzero element in the kernel 
gives a non trivial relation, over the field $\Q(\zeta_{\prod_{t\not=r}  q_{t}})$, 
between the elements of a $\Q$-basis of 
$\sum_{j=1}^{\ell_t}(i\pi)^{-k_r}V_{\Q,k_r}^+(q_{r,j})$. But this contradicts 
the linear disjointness of the fields $\Q(\zeta_{\prod_{t\not=r}  q_{t} } )$ and 
$\Q(\zeta_{q_{r}})$. Thus, the dimension of the $\Q$ vector space 
on the right is the product 
$\prod_{t=1}^r\dim_{\Q}\left(\sum_{j=1}^{\ell_t}(i\pi)^{-k_t}V_{\Q,k_t}^+(q_{t, j})\right)$ 
which is equal to $2^{-r}\prod_{t=1}^{r} \left(\sum_{j=1}^{\ell_t}\varphi(q_{t, j}) - (\ell_t - 1)(1+(-1)^{k_t})\right)$ 
by Theorem \ref{cor1}.
\qed

\medskip

\section{Appendix: alternative proof of Proposition \ref{prop1}.}
\subsection{The case when $k$ is odd}
\begin{proof}
Set $q= q_1 \cdots q_{\ell}$. Thus, applying \eqref{dcot} and \eqref{ecot}, we have
\begin{eqnarray*}
&&
-i^{2k+1} (\zeta_{q} - \zeta_{q}^{-1})\cot^{(2k)} (\pi a_j/ q_j)  \\
&=&
(\zeta_{q} - \zeta_{q}^{-1}) 
\sum_{a, b \ge 0 \atop a + 2b = 2k +1} 
 \beta_{a, b}^{(2k+1)} (-i)^a (-1)^{b} (\cot \frac{\pi a_j}{q_j})^a
~(1 + (\cot  \frac{\pi a_j}{q_j})^2)^b \\
&=& 
(\zeta_{q} - \zeta_{q}^{-1})
\sum_{a, b \ge 0 \atop a + 2b = 2k +1} 
 (-1)^b \beta_{a, b}^{(2k+1)} \left( -i \cot \frac{\pi a_j}{q_j} \right)^a
~\left(1 - (-i \cot  \frac{\pi a_j}{q_j})^2 \right)^b \\
&=&
(\zeta_{q} - \zeta_{q}^{-1})   
\sum_{a, b \ge 0 \atop a + 2b = 2k +1 } 
(-1)^b\beta_{a, b}^{(2k+1)} \left( \frac{\zeta^{a_j}_{q_j} + 1}{\zeta^{a_j}_{q_j} - 1}  \right)^a
~\left(1 - (\frac{\zeta^{a_j}_{q_j} + 1}{\zeta^{a_j}_{q_j} - 1} )^2 \right)^b. 
\end{eqnarray*}
This implies that 
$i(\zeta_{q} - \zeta_{q}^{-1})\cot^{(2k)} (\frac{\pi a_{j}}{q_{j}}) \in \Q(\zeta_{q})^+$ 
for $a_j \in T_{q_j}, ~1 \le j \le \ell$.
The proposition then follows for $\ell =1$ along the lines of the proof of Okada and  Murty-Saradha. 
Suppose that the proposition is true for any natural number $1 \leq n < \ell$.
We want to show that 
\begin{equation}\label{der-cot}
\bigcup_{ 1 \leq j \leq \ell} \{ \cot^{(2k)} (\frac{\pi a_j}{q_j}) ~ : ~ a_j \in T_{q_j}\}
\end{equation}
is linearly independent over $\K$. By the given hypothesis, we have 
$\K \cap \Q(\zeta_{q_1\cdots q_{\ell} })^+ = \Q$.
Appealing to \thmref{kf}, it is now sufficient to show that the numbers in the set \eqref{der-cot}
are $\Q$ linearly independent. There exist rational numbers 
$\alpha_{a_{j}}$ for $a_j \in T_{q_j}, 1 \leq j \leq \ell$ such that 
$$
\sum_{ 1 \leq j \leq \ell} \sum_{a_j \in T_{q_j}} \alpha_{a_{j}} \cot^{(2k)} (\frac{\pi a_j}{q_j}) 
~=~ 0.
$$
This implies that
$$
\sum_{ 1 \leq j < \ell} \sum_{a_j \in T_{q_j}} \alpha_{a_{j}} \cot^{(2k)} (\frac{\pi a_j}{q_j}) 
~~=~~
 - \sum_{a_{\ell} \in T_{q_\ell}} \alpha_{a_{\ell}} \cot^{(2k)} (\frac{\pi a_{\ell}}{q_{\ell}}).
$$
Alternative to the identity \eqref{dcot}, one can write
$$
\cot^{(h-1)}(z) = i^h \left(2X\frac{d}{dX} \right)^{h-1}\left( \frac{X+1}{X-1}\right)|_{X = e^{2iz}}
$$ 
for all $h$  since $\cot(z) = i \frac{e^{2iz}+1}{e^{2iz}-1}$.  Evaluating the above expression
at $z=\pi a_j/q_j$, one gets $i^h \cot^{(h-1)}(\pi a_j/q_j) \in \Q(\zeta_{q_j})\cap i^h\R$.
Since $h=2k+1$ is an odd integer, it then follows that
$$
i \sum_{ 1 \leq j < \ell} \sum_{a_j \in T_{q_j}} \alpha_{a_{j}} \cot^{(2k)} (\frac{\pi a_j}{q_j}) 
~=~ 
- i \sum_{a_{\ell} \in T_{q_\ell}} \alpha_{a_{\ell}} \cot^{(2k)} (\frac{\pi a_{\ell}}{q_{\ell}})
~\in~ 
\Q(\zeta_{q_1 \cdots q_{\ell - 1}}) \cap \Q(\zeta_{q_{\ell}})
~=~ \Q.
$$
Since a purely imaginary number is a rational number if and only if it is $0$, we have
\begin{align*}
\begin{split}
\sum_{ 1 \leq j < \ell} \sum_{a_j \in T_{q_j}} \alpha_{a_{j}}  \cot^{(2k)} (\frac{\pi a_j}{q_j})
~=~ 
- \sum_{a_{\ell} \in T_{q_\ell}} \alpha_{a_{\ell}} \cot^{(2k)} (\frac{\pi a_{\ell}}{q_{\ell}}) 
&~=~ 0.
\end{split}
\end{align*}
Applying induction hypothesis, we get that $\alpha_{a_{j}} = 0$ 
for all $a_j \in T_{q_j},  1 \leq j \leq \ell$. 
This completes the proof of the proposition.
\end{proof}

\subsection{The case when $k$ is even}
\begin{proof}
Applying \eqref{dcot}, we see that
$\cot^{(2k-1)} (\frac{\pi a_{j}}{q_{j}}) \in \Q(\zeta_{q_j})^+$ 
for $a_j \in T_{q_j}, ~1 \le j \le \ell$. The proposition then follows for $\ell =1$ from the 
works of Okada  and  Murty-Saradha. 
Suppose that the proposition is true for any $1\leq n < \ell$.  We want to show that
the numbers
\begin{equation}\label{d-co}
\{ \cot^{(2k-1)} (\frac{\pi a_{1}}{q_{1}}) ~ : ~ a_{1} \in T_{q_1}\} 
\bigcup_{ 1 < j  \leq  \ell} \{ \cot^{(2k-1)} (\frac{\pi a_j}{q_j}) ~ : ~ a_j \in S_{q_j} \} 
\end{equation}
are linearly independent over $\K$. Since $\K \cap \Q(\zeta_{q_1\cdots q_{\ell}})^+ = \Q$,
applying \thmref{kf}, we see that it is sufficient to prove that
the numbers in the set \eqref{d-co} are $\Q$ linearly independent. 
There exist rational numbers $\alpha_{a_{j}} $ for 
$a_1 \in T_{q_1}, a_j \in S_{q_j}, 1 < j \leq \ell$ 
such that 
$$
\sum_{a_1 \in T_{q_1}} \alpha_{a_{1}} \cot^{(2k-1)} (\frac{\pi a_1}{q_1}) 
~+~ 	
\sum_{ 1 < j \leq \ell} \sum_{a_j \in S_{q_j}} \alpha_{a_{j}} \cot^{(2k-1)} (\frac{\pi a_j}{q_j}) 
~=~ 0.
$$
Then
\begin{eqnarray*}
&&
\sum_{a_1 \in T_{q_1}} \alpha_{a_{1}} \cot^{(2k-1)} (\frac{\pi a_1}{q_1}) 
~~+~~ 
\sum_{ 1 < j < \ell} \sum_{a_j \in S_{q_j}} \alpha_{a_{j}} \cot^{(2k-1)} (\frac{\pi a_j}{q_j}) \\
&=& 
 \sum_{a_{\ell} \in S_{q_\ell}}  \alpha_{a_{\ell}} \cot^{(2k-1)} (\frac{\pi a_{\ell}}{q_{\ell}})  
~\in~ 
\Q(\zeta_{q_1 \cdots q_{\ell-1}})  \cap \Q(\zeta_{q_{\ell}})
~=~ \Q.
\end{eqnarray*}
Let us call this rational number $\beta$. If $\beta \ne 0$, then without loss of generality, 
we may assume that $\beta = - 1$. So
\begin{align}\begin{split}\label{stu}
 \sum_{a_{\ell} \in S_{q_\ell}} \alpha_{a_{\ell}} \cot^{(2k-1)} (\frac{\pi a_{\ell}}{q_{\ell}}) 
~=~ 1.
\end{split}
\end{align}
Let us denote $\Q(\zeta_{j})$ by $\F_j$. Let $\T_{\F_j/\Q}(\alpha)$ denotes
the trace of $\alpha$ over $\Q$ for $\alpha \in \F_j$. For any $a_{j} \in T_{j}$, 
using identities \eqref{dcot} and \eqref{ecot}, we get
\begin{eqnarray}\label{tr}
tr_j 
&=& 
\T_{\F_j /\Q} \left( \cot^{(2k-1)} (\frac{\pi a_{j}}{q_{j}})\right) \nonumber \\
&=&
\sum_{ 1 \leq a_{j} \leq q_{j} \atop (a_{j},q_{j}) = 1  }\cot^{(2k-1)} (\frac{\pi a_{j}}{q_{j}})
~=~ 
2 \sum_{a_{j} \in T_{q_j}} \cot^{(2k-1)} (\frac{\pi a_{j}}{q_{j}}).
\end{eqnarray}
It now follows from \eqref{stu} and \eqref{tr} that
\begin{align*}
\begin{split}
 \sum_{a_{\ell} \in S_{q_\ell}} \alpha_{a_{\ell}} 
\cot^{(2k-1)} (\frac{\pi a_{\ell}}{q_{\ell}}) 
= 
\frac{2}{tr_{\ell}}\left( \sum_{a_{\ell} \in T_{q_\ell}} \cot^{(2k-1)} 
(\frac{\pi a_{\ell}}{q_{\ell}})\right).
\end{split}
\end{align*}
Hence we have
\begin{align*}
\begin{split}
\frac{-2}{tr_{\ell}}\cot^{(2k-1)} (\frac{\pi a_{1}}{q_{1}}) 
+  \sum_{a_{\ell} \in S_{q_\ell}} (\alpha_{a_{\ell}} 
- \frac{2}{tr_{\ell}} )\cot^{(2k-1)} (\frac{\pi a_{\ell}}{q_{\ell}}) = 0,
\end{split}
\end{align*}
a contradiction to \thmref{okada}. This implies that $\beta = 0$.
Applying induction hypothesis, we then obtain 
$\alpha_{a_{j}} = 0 $ for $a_1 \in T_{q_1}, ~a_j \in S_{q_j}, 1 < j \leq \ell$.
\end{proof}

\bigskip

\noindent
{\bf Acknowledgments.} 
Research of this article was partially supported by Indo-French 
Program in Mathematics (IFPM).  All authors would like to thank IFPM 
for financial support. The first and the second author also thank DAE
number theory plan project. The second author would also like to thank 
Project ANR JINVARIANT for partial financial support.


\begin{thebibliography}{100}
\bibitem{BBW} 
A. Baker,  B.J. Birch and  E. A Wirsing,
{\em On a problem of Chowla},
J. Number Theory {\bf 5} (1973), 224--236. 


\bibitem{BR}
K. Ball and T. Rivoal,
{\em Irrationalité d'une infinité de valeurs de la fonction zêta 
aux entiers impairs}, 
Invent. Math. {\bf 146} (2001), 193--207. 


\bibitem{chowla}
 S. Chowla,
 {\em The nonexistence of nontrivial linear relations between the roots of 
 a certain irreducible equation}, 
 J. Number Theory {\bf 2} (1970), no. 2, 120--123.
 
 
 \bibitem{CC} 
P. Chowla and S. Chowla, {\em On irrational numbers}, 
 Skr. K. Nor. Vidensk. Selsk. (Trondheim) 3, (1982), 1-5. 
(See also S. Chowla, Collected Papers, Vol 3,
pp. 1383-1387, CRM, Montreal, 1999.)


\bibitem{MC}
P.  M. Cohn, 
Algebra, Second edition, vol. 3, {\em John Wiley \& Sons}, 1991.


\bibitem{SF}
S.  Fischler, 
{\em Irrationality of values of L-functions of Dirichlet characters},
J. Lond. Math. Soc. (2) {\bf 101} (2020), no. 2, 857--876.
 

\bibitem{FSZ}
S. Fischler, J. Sprang and W. Zudilin, 
{\em Many odd zeta values are irrational}, 
Compos. Math. {\bf 155} (2019), no. 5, 938--952.  
 
 
 \bibitem{KG}  
K. Girstmair, 
{\em Letter to the editor}, 
J. Number Theory {\bf 23} (1986), p. 405.


\bibitem{GMR}
S. Gun, M. Ram Murty and P. Rath,
{\em On a conjecture of Chowla and Milnor},
Canad. J. Math. {\bf 63} (2011), 1328--1344.


\bibitem{GK}
S. Gun and N. Kandhil,
{\em On an extension of a question of Baker}, 
Int. J. Number Theory, \\
https://doi.org/10.1142/S1793042123500173.


\bibitem{YH}
Y. Hamahata, 
{\em Okada's theorem and multiple Dirichlet series}, 
Kyushu J. Math {\bf 74} (2020), 429--439.


\bibitem{HH}
H. Hasse, 
{\em On a question of S. Chowla},
Acta Arithmetica {\bf 18} (1971), 275--280.


\bibitem{JH}
H. Jager and H.W. Lenstra,
{\em Linear independence of cosecant values}, 
Nieuw Arch. Wisk {\bf 23} (1975), no. 3, 131--144.


\bibitem{NR} 
N. Kandhil and P. Rath, 
{\em Around a question of Baker}, to appear in Proc.
of Subbarao Symposium to be published by the Fields institute. 


\bibitem{LY}
L. Lai and P. Yu, 
{\em  A note on the number of irrational odd zeta values},
Compos. Math. {\bf 156} (2020), no. 8, 1699--1717.


\bibitem{MS1}  
M. Ram Murty and N. Saradha,
{\em Special values of the polygamma functions,} 
Int. J. Number Theory {\bf 5} (2009), no. 2, 257--270. 


\bibitem{MM1}  
M. Ram Murty and V. Kumar Murty,
{\em A problem of Chowla revisited,} 
J. Number Theory {\bf 131} (2011), no. 9, 1723--1733.


\bibitem{MIL}
J. Milnor, 
{\em On polylogarithms, Hurwitz zeta functions,
and their Kubert identities}, Enseignement Math. (2) 29 (1983),
no. 3-4, 281--322. 


\bibitem{JN}
J. Neukirch,
{\em Algebraic Number Theory}, Springer, 1999. 


\bibitem{TO} 
T. Okada, 
{\em On an extension of a theorem of S. Chowla,}  
 Acta Arithmetica {\bf 38} (1981), no. 4, page 341--345.


\bibitem{RZ}
T. Rivoal and W. Zudilin, 
{\em Diophantine properties of numbers related to 
Catalan's constant},
 Math. Ann. {\bf 326} (2003), no. 4, 705--721.


\bibitem{KW}
K. Wang, 
{\em On a theorem of S. Chowla},
J. Number Theory {\bf 15} (1982), 1--4.


\end{thebibliography}
\end{document}